\documentclass[10pt,onesite,draft]{article} 
\usepackage{amsmath,amsxtra,amssymb,latexsym, amsfonts, indentfirst}
\usepackage[mathscr]{eucal}
\newfam\cyrfam

\def\ovl{\overline}

\newenvironment{proof}{%
\par\addvspace{6pt plus3pt minus2pt}%
\noindent{\bfseries\itshape\textit Proof.}\ignorespaces} {%
\if@halmos\halmos\fi
\par\addvspace{6pt plus3pt minus2pt} }
\begin{document}
\parskip 4pt
\large
\setlength{\baselineskip}{15 truept}
\setlength{\oddsidemargin} {0.5in}
\overfullrule=0mm
\def\bfh{\vhtimeb}
\date{}
\title{\bf \large  FILTER-REGULAR SEQUENCES  AND MIXED  MULTIPLICITIES$^1$}
\def\b{\vntime}
\author {{\small  Nguyen Tien Manh  and  Duong Quoc Viet}\\    
{\small Department of Mathematics, Hanoi  University of Education}\\
{\small 136-Xuan Thuy Street, Hanoi, Vietnam}\\
{\small\it E-mail: duongquocviet@fmail.vnn.vn}}
\maketitle
\centerline{
\parbox[c]{12 cm}{
\small  ABSTRACT:   Let $S$ be a finitely generated standard multigraded algebra over an Artinian local ring $A$; $M$   a finitely generated multigraded $S$-module.  This paper answers to the question  when   mixed multiplicities of $M$ are positive and characterizes them in terms of  lengths of $A$-modules. As an application, we get interesting results on mixed multiplicities of ideals, and recover some early results in [Te] and [TV]. }}

\vspace{24pt}
\centerline{\Large\bf 1. Introduction}$\\$
Throughout this paper, let $(A,\frak{m})$ denote an Artinian local ring with  maximal ideal $\frak{m}$, infinite residue field $k = A/\frak{m}$; $S=\bigoplus_{n_1,\ldots,n_d\ge 0}S_{(n_1,\ldots,n_d)}$ $(d > 0)$ a finitely generated standard $d$-graded algebra over $A$ (i.e., $S$ is generated over $A$ by elements of total degree 1); $M=\bigoplus_{n_1,\ldots,n_d\ge 0}M_{(n_1,\ldots,n_d)}$   a finitely generated $d$-graded $S$-module.
Set  
$$\begin{array}{lll}&{\frak a}:
{\frak b}^\infty&=\bigcup_{n\ge 0}(\frak a:{\frak b}^n);\\&S^\triangle&=\bigoplus_{n\ge 0}S_{(n,\ldots,n)};\;S^\triangle_+=\bigoplus_{n> 0}S_{(n,\ldots,n)};\\&S_+&=\bigoplus _{n_1+\cdots+n_d> 0}S_{(n_1,\ldots,n_d)};\;S_{++}=\bigoplus_{n_1,\ldots,n_d> 0}S_{(n_1,\ldots,n_d)};\\&S_i&=S_{(0,\ldots,{\underbrace{1}_i},\ldots,0)}\;(i=1,\ldots,d);\\&M^\triangle&=\bigoplus_{n\ge 0}M_{(n,\ldots,n)};\;\ell=\dim M^\triangle.\end{array}$$
\footnotetext[1]{\begin{itemize}
\item[ ]{\it Mathematics  Subject  Classification} (2000): Primary 13H15. Secondary 13A02, 13E05, 13E10, 14C17. 
\item[ ]$ Key\; words \;  and \; phrases:$ Artinian  ring,  multiplicity,  multigraded  module, filter-regular  sequence.
\end{itemize}}Denote by $\text{Proj}\; S$ the set of the homogeneous prime ideals of $S$ which do not contain $S_{++}$. 
Set $$\text{Supp}_{++}M=\{P\in \text{Proj}\; S\;|\;M_P\ne 0\}.$$   By [HHRT, Theorem 4.1] and Remark 3.1,  $\dim \text{Supp}_{++}M=\ell-1$ and $l_A[M_{(n_1,\ldots,n_d)}]$ is a polynomial of degree $\ell-1$ for all large $n_1,\ldots,n_d$ 
(see Remark 3.1, Section 3).  The terms of total degree $\ell-1$ in this polynomial have the form 

\vspace*{6pt}

$$\sum_{k_1\:+\:\cdots\:+\:k_d\;=\;\ell-1}e(M;k_1,\ldots,k_d)\dfrac{n_1^{k_1}\cdots n_d^{k_d}}{k_1!\cdots k_d!}.$$ 

\vspace*{6pt}
\hspace*{-0.5cm}Then $e(M;k_1,\ldots,k_d)$ are   non-negative integers not all zero, called the {\it  mixed multiplicity of type $(k_1,\ldots,k_d)$ of $M$} [HHRT]. 

\vspace*{10pt}

In the case that $(R, \frak n)$  is  a  Noetherian   local ring with  maximal ideal $\frak{n},$ $J$ is an $\frak n$-primary ideal, $I_1,\ldots, I_s$ are ideals of $R.$ Then 

\vspace*{6pt}

$$ T  =\bigoplus_{n, n_1,\ldots,n_s\ge 0}\dfrac{J^nI_1^{n_1}\cdots I_s^{n_s}}{J^{n+1}I_1^{n_1}\cdots I_s^{n_s}}$$

\vspace*{6pt}
\hspace*{-0.5cm}is a finitely generated standard multigraded algebra over Artinian local ring  $R/J.$ Mixed multiplicities of $T$ are called 
mixed multiplicities of ideals $J,I_1,\ldots, I_s$ (see [Ve] or [HHRT]).

\vspace*{6pt}
The theory mixed multiplicities of $\frak{n}$-primary ideals was introduced   by Risler and Teissier in 1973 [Te], by Rees in 1984 [Re].  
 In past years, the positivity and the
relationship between  mixed multiplicities  and  Hilbert-Samuel multiplicity of ideals have attracted much attention (see e.g. [Sw], [Ro], [KR1, KR2], [KV], [Vi], [Tr1], [Tr2], [MV], [TV]).

\vspace*{10pt}

We turn now to the positivity of mixed multiplicities of a multigraded module. In the geometric context, this problem  was introduced by Kleiman and Thorup in 1996 [KT].  In the algebraic setting, it was investigated by Trung in 2001 [Tr2] in the case of a bigraded ring.

\vspace*{10pt}

Under a totally algebraic point of view, 
this paper answers to the question  when  mixed multiplicities of multigraded modules are positive and characterizes these  mixed multiplicities in terms of lengths of modules. Our approach is based on the properties of   filter-regular sequences. The notion of filter-regular sequences was introduced by Stuckrad and Vogel in their book [SV]. The theory of filter-regular sequences became an important tool to study some classes of singular rings and has been continually developed (see e.g. [Tr1], [BS], [Hy],[Tr2]). As one might expect, we obtain the following result.\\\\
{\bf Main theorem} (Theorem 3.4){\bf.} {\it Let $S$ be a finitely generated standard $d$-graded algebra over  an Artinian local ring $A$ and $M$  a finitely generated $d$-graded $S$-module such that $S_{(1,1,\ldots,1)}$ is not contained in $ \sqrt{\mathrm{Ann}_S M}$. Set $\ell=\dim M^\triangle.$ Then the following statements hold.
\begin{itemize}
\item[$\mathrm{(i)}$] $e(M;k_1,\ldots,k_d)\ne 0 $\; if and only if  there exists an $S_{++}$-filter-regular sequence $x_1,\ldots,x_{\ell-1}$  with respect to  $M$ consisting of $k_1$ elements of $S_1,\ldots,k_d$ elements of  $S_d$ and 
$\dim \ovl M^{\triangle} =1$, where $\ovl M=\dfrac{M}{(x_1,\ldots,x_{\ell-1})M}$.
\item[$\mathrm{(ii)}$] Suppose that $e(M;k_1,\ldots,k_d)\ne 0 $ and $x_1,\ldots,x_{\ell-1} $ is  an $S_{++}$-filter-regular sequence  with respect to  $M$ consisting of $k_1$ elements of $S_1,\ldots,k_d$ elements of  $S_d$. Set $\ovl M=\dfrac{M}{(x_1,\ldots,x_{\ell-1})M}\; and\; r=r(S^{\triangle}_+;\ovl M^{\triangle})$ the reduction number of $S^{\triangle}_+$ with respect to  $\ovl M^{\triangle}.$  Then $ e(M;k_1,\ldots,k_d)=l_A\biggl[\dfrac{\ovl M^\triangle_n}{(0_{\ovl M^\triangle}:S_+^{\triangle\infty})_n}\biggl]\;\text{for\; all\;} n\ge r.$
\end{itemize}}
 As interesting consequences of main result, we get results on mixed multiplicities of ideals (Proposition 4.4 and Theorem 4.5), and  recover the results of Risler and Teissier [Te](see Remark 4.8),  Trung and Verma [TV](see Remark 4.6).   

 \vspace*{12pt}
\centerline{\Large\bf2. On Filter-Regular Sequences }

\vspace*{12pt}

Recently, filter-regular sequences have been used to investigate mixed multiplicities and joint reduction numbers in bigraded rings (see e.g. [Hy], [Tr2]). This section gives some properties of filter-regular sequences in  multigraded modules.\\\\
{\bf Definition 2.1.} {\it Let $S=\bigoplus_{n_1,\ldots,n_d\ge 0}S_{(n_1,\ldots,n_d)}$  be a finitely generated standard $d$-graded algebra  over an Artinian local ring $A$; $M=\bigoplus_{n_1,\ldots,n_d\ge 0}M_{(n_1,\ldots,n_d)}$   a finitely generated $d$-graded $S$-module such that $S_{(1,1,\ldots,1)}$ is not contained in $ \sqrt{\mathrm{Ann}_SM}$. A homogeneous element $x\in S$ is called an $S_{++}$-filter-regular element with respect to  $M$  if
$x\notin P$ for any $P\in \mathrm{Ass}_SM$ and $P$ does not contain   $S_{++}.$ That means 

\centerline {$x\notin\bigcup_{S_{++}\nsubseteq P,\; P\in \mathrm{Ass}_SM}P.$}
\hspace*{-0.6cm}Let  $x_1,\ldots, x_t$ be homogeneous elements in $S$. We call that $x_1,\ldots, x_t$ is an $S_{++}$-filter-regular sequence  with respect to  $M$ if 
 $x_i$ is an $S_{++}$-filter-regular element with respect to $\dfrac{M}{(x_1,\ldots, x_{i-1})M}$ for all $i = 1,\ldots, t.$}\\

Now, we  briefly  give some comments  on filter-regular sequences  of a finitely generated  multigraded  module.\\

\hspace*{-0.6cm}{\bf Note:} Since $\mathrm{Ass}_S[M/(0_M:S_{++}^{\infty})]= \{P\in \mathrm{Ass}_S M\;|\;S_{++}\nsubseteq P\},$ it follows that a homogeneous element $x\in S$ is  an $S_{++}$-filter-regular element if and only if $x$ is  non-zero-divisor in $M/(0_M:S_{++}^{\infty}).$ Moreover, we have the following notes. 
\begin{itemize}
\item[$\mathrm{(i)}$] It is easy to see that a homogeneous element $x\in S$ is  an $S_{++}$-filter-regular element with respect to $M$ if and only if $0_M:x\subseteq 0_M:S_{++}^{\infty}.$  

\item[$\mathrm{(ii)}$] If  $S_{(1,1,\ldots,1)}\subseteq \sqrt{\mathrm{Ann}_S M}$ then $0_M:S_{++}^{\infty} = M.$  Hence for any homogeneous element $x\in S$, we always have  $0_M:x \subseteq 0_M:S_{++}^{\infty}.$  This only obstructs and does not carry usefulness. That is why in Definition 2.1, one has to exclude the case that $S_{(1,1,\ldots,1)}\subseteq \sqrt{\mathrm{Ann}_S M}.$\\

\item[$\mathrm{(iii)}$] An $S_{++}$-filter-regular sequence $x_1,\ldots, x_t$ with respect to $M$ is a maximal $S_{++}$-filter-regular sequence if  $$S_{(1,1,\ldots,1)}\nsubseteq \sqrt{\mathrm{Ann}_S [M/(x_1,\ldots, x_{t-1})M]}$$  and  $ S_{(1,1,\ldots,1)}\subseteq \sqrt{\mathrm{Ann}_S [M/(x_1,\ldots, x_t)M]}.$\\
\end{itemize}
 
Let $\bf{m}$ $=(m_1,\ldots,m_d),\bf{n}$ $=(n_1,\ldots,n_d)\in \mathbb{Z}^d$. We recall that $\bf m\le \bf n$ if 
 $m_i\le n_i$ for all $i=1,\ldots,d$. 
\\\\
{\bf Lemma 2.2.} {\it Let $N$ be an $d$-graded $S$-submodule of $M$. Then there exist  positive integers $u_1,\ldots,u_d$  such that $N_{(n_1,\ldots,n_d)}=S_{(n_1-u_1,\ldots,n_d-u_d)}N_{(u_1,\ldots,u_d)}$ for all $n_1\ge u_1,\ldots,n_d\ge u_d$.}

\begin{proof}\; Since $M$ is a finitely generated $S$-module and $S$ is a Noetherian ring, $N$ is also a finitely generated $S$-module. It implies that there exist positive integers $u_1,\ldots,u_d$ such that $N$ is generated by 
$\bigcup_{v_1\le u_1,\ldots,v_d\le u_d}N_{(v_1,\ldots,v_d)}.$  Set $\bf u$ $=(u_1,\ldots,u_d)$, $\bf v$ $=(v_1,\ldots,v_d)$, $\bf n$ =$(n_1,\ldots,n_d)$. We have $$N_{\bf n}=\sum_{{\bf v}\le {\bf u}}S_{\bf n-v}N_{\bf v}\supseteq S_{\bf n-u}N_{\bf u}$$ for all ${\bf n}\ge {\bf u}$. Since $S$ is a  finitely generated standard $d$-graded ring,  $S_{\bf m}S_{\bf n}=S_{\bf m+n}$ for all $\bf m,n.$  It follows that   $S_{\bf n-v}N_{\bf v}= S_{\bf n-u}(S_{\bf u-v}N_{\bf v})\subseteq S_{\bf n-u}N_{\bf u}$ for all $\bf v\le \bf u$. Hence $N_{\bf n}\subseteq S_{\bf n-u}N_{\bf u}$  for all ${\bf n}\ge {\bf u}$.  So  $N_{\bf n}= S_{\bf n-u}N_{\bf u}$ for all ${\bf n}\ge {\bf u}$. $\blacksquare$
\end{proof}

\hspace*{-0.6cm}{\bf Remark 2.3.} Since $M$ is Noetherian, there exists a positive integer $n$ such that $0_M:S_{++}^{\infty}=0_M:S_{++}^{n}.$    
By Lemma 2.2, there exist  positive integers $u_1,\ldots,u_d$  such that $$[0_M:S_{++}^{n}]_{(n_1,\ldots,n_d)}=S_{(n_1-u_1,\ldots,n_d-u_d)}[0_M:S_{++}^{n}]_{(u_1,\ldots,u_d)}$$ for all $n_1\ge u_1,\ldots,n_d\ge u_d$.
This fact follows that  $[0_M:S_{++}^{n}]_{(n_1,\ldots,n_d)}=0$ for all $n_1\ge n+u_1,\ldots,n_d\ge n+u_d$. Hence $(0_M:S_{++}^{\infty})_{(n_1,\ldots,n_d)}=0$ for all large $n_1,\ldots,n_d$.\\

The following characteristic of  filter-regular sequences will often be used.\\\\
{\bf Proposition 2.4.} {\it Let $x\in S$ be a homogeneous element. Then $x$ is an $S_{++}$-filter-regular element with respect to  $M$ if and only if $(0_M:x)_{(n_1,\ldots,n_d)}=0$ for  all large $n_1,\ldots,n_d.$}

\begin{proof}\: Suppose  that $x$ is an $S_{++}$-filter-regular element with respect to $M$. 
Then $0_M:x\subseteq 0_M:S_{++}^{\infty}$ by Note (i). 
By Remark 2.3, $$(0_M:x)_{(n_1,\ldots,n_d)}\subseteq (0_M:S_{++}^{\infty})_{(n_1,\ldots,n_d)}=0$$ for all large $n_1,\ldots,n_d$.
Now assume that $(0_M:x)_{(n_1,\ldots,n_d)}=0$ for  all large $n_1,\ldots,n_d.$
It is enough to prove that 
$(0_M:x)_{(v_1,\ldots,v_d)}\subseteq0_M:S_{++}^{\infty}$ for all $v_1,\ldots,v_d.$ 
Indeed, if $a\in (0_M:x)_{(v_1,\ldots,v_d)}$ then $xS_{(n_1,\ldots,n_d)}a=0$ for every $n_1,\ldots,n_d$. Hence $$S_{(n_1,\ldots,n_d)}a\subseteq (0_M:x)_{(n_1+v_1,\ldots,n_d+v_d)}=0$$ for large $n_1,\ldots,n_d$. It implies that  $a\in 0_M:S_{++}^{n}$ for all large $n.$ Hence $a\in 0_M:S_{++}^{\infty}.$ So 
$(0_M:x)_{(v_1,\ldots,v_d)}\subseteq0_M:S_{++}^{\infty}$ for all $v_1,\ldots,v_d.$ 
 Hence $x$ is an $S_{++}$-filter-regular element with respect to  $M$.
$\blacksquare$
\end{proof}

\hspace*{-0.6cm}{\bf Remark 2.5.} 
 Suppose that $x\in S_{(a_1,\ldots,a_d)}$ is  an $S_{++}$-filter-regular element with respect to $M$. Consider \\\centerline{$\lambda_x:M_{(n_1-a_1,\ldots,n_d-a_d)}\longrightarrow xM_{(n_1-a_1,\ldots,n_d-a_d)}, y\mapsto xy.$} It is clear that $\lambda_x$ is  surjective and $\mathrm{Ker}\lambda_x =(0_M:x)_{(n_1-a_1,\ldots,n_d-a_d)}$. Since $x$ is an $S_{++}$-filter-regular element with respect to $M$, $\mathrm{Ker}\lambda_x =0$ for all large $n_1,\ldots,n_d$ by Proposition 2.4. Therefore, $xM_{(n_1-a_1,\ldots,n_d-a_d)}\cong M_{(n_1-a_1,\ldots,n_d-a_d)}$    and   $$\begin{array}{lll}l_A[(M/xM)_{(n_1,\ldots,n_d)}]&=&l_A[M_{(n_1,\ldots,n_d)}]-l_A[xM_{(n_1-a_1,\ldots,n_d-a_d)}]\\&=&l_A[M_{(n_1,\ldots,n_d)}]-l_A[M_{(n_1-a_1,\ldots,n_d-a_d)}]\end{array}$$ for all large $n_1,\ldots,n_d$.\\

Next, we show that  the existence  of filter-regular sequences is universal.\\\\
{\bf Proposition 2.6.} {\it Assume that $S_{(1,1,\ldots,1)}$ is not contained in $ \sqrt{\mathrm{Ann}_S M}$. 
 Then for each $i=1,\ldots,d$, there exists an $S_{++}$-filter-regular element $x\in S_i$ with respect to  $M$.}

\begin{proof}\: Set $M^*=M/(0_M:S_{++}^{\infty}).$  Note that
$$\mathrm{Ass}_S M^*=\{P\in \mathrm{Ass}_S M\;|\;S_1\cdots S_d = S_{(1,1,\ldots,1)}\nsubseteq P\},$$ $S_i \nsubseteq P \in \mathrm{Ass}_S M^*$ for all $i= 1,\ldots,d.$ Since $k$ is an infinite field and $\mathrm{Ass}_S M^*$ is finite. Hence for each $i=1,\ldots,d,$  there exists   $x\in S_i\setminus\bigcup_{P\in\mathrm{Ass}_S M^* }P.$ On the other hand $S_{(1,1,\ldots,1)}\nsubseteq \sqrt{\mathrm{Ann}_S M}$, $x$ is an $S_{++}$-filter-regular element with respect to $M$.$\blacksquare$
\end{proof}
Suppose that $S_{(1,1,\ldots,1)}$ is not contained in $ \sqrt{\mathrm{Ann}_S M}$. Then  $S_{(u,\ldots,u)}M\ne 0$ for all $u > 0$. Hence for any $u >0$, there exist $m_1,\ldots,m_d$ such that $S_{(u,\ldots,u)}M_{(m_1,\ldots,m_d)}\ne 0.$ Since $S_{(u,\ldots,u)}M_{(m_1,\ldots,m_d)}\subseteq M_{(m_1+u,\ldots,m_d+u)},$ $M_{(m_1+u,\ldots,m_d+u)}\ne 0.$ So for any $u > 0$, there exist $n_1,\ldots,n_d \ge u$ such that $M_{(n_1,\ldots,n_d)}\ne 0.$ 
Now, assume that for each $u >0$, there exist  $n_1,\ldots,n_d\ge u$ such that $M_{(n_1,\ldots,n_d)}\ne 0.$ We will show that $S_{(1,1,\ldots,1)}\nsubseteq \sqrt{\mathrm{Ann}_S M}.$ 
Indeed, if $S_{(1,1,\ldots,1)}\subseteq \sqrt{\mathrm{Ann}_S M}$ then $S_{(n,\ldots,n)}M=0$ for some $n$. By Lemma 2.2, there exist  positive integers $u_1,\ldots,u_d$  such that \\\centerline{$M_{(n_1,\ldots,n_d)}=S_{(n_1-u_1,\ldots,n_d-u_d)}M_{(u_1,\ldots,u_d)}$} for all $n_1\ge u_1,\ldots,n_d\ge u_d$. Hence \\\centerline{$M_{(n_1,\ldots,n_d)}=S_{(n_1-n-u_1,\ldots,n_d-n-u_d)}S_{(n,\ldots,n)}M_{(u_1,\ldots,u_d)}=0$} for all $n_1\ge n+u_1,\ldots,n_d\ge n+u_d$. This contracdicts the hypothesis. So $S_{(1,1,\ldots,1)}\nsubseteq \sqrt{\mathrm{Ann}_S M}.$ From these facts, we  get the following proposition.\\\\ 
{\bf Proposition 2.7.} {\it Let $S$  be a finitely generated standard $d$-graded algebra over an Artinian local ring $A$; $M$   a finitely generated $d$-graded $S$-module. Then the following conditions are equivalent:
\begin{itemize}
\item[$\mathrm{(i)}$] $S_{(1,1,\ldots,1)}$ is  contained in $ \sqrt{\mathrm{Ann}_S M}$.
\item[$\mathrm{(ii)}$] $M_{(n_1,\ldots,n_d)}= 0$ for all large $n_1,\ldots,n_d.$
\end{itemize}}

Let $S=\bigoplus_{n\ge 0}S_n$  be a finitely generated standard graded algebra over an Artinian local ring $A$ and  $J$  a homogeneous ideal of $S$ is generated by elements of degree 1. Let $M=\bigoplus_{n\ge 0}M_n$ be a finitely generated graded $S$-module. Set $S_+=\bigoplus_{n> 0}S_n$.  We call $J$  {\it a reduction of $S_+$ with respect to $M$} if $(JM)_n=M_n$ for all large  $n$.  The least integer $n$ such that $(JM)_{n+1}=M_{n+1}$ is called {\it the reduction number of $S_+$ with respect to $J$ and $M$} [NR]. We denote this integer by  $r_J(S_+;M).$ A  reduction $J$ of $S_+$ with respect to  $M$ is called  {\it a minimal reduction } if it does not properly contain  any other reduction of $S_+$ with respect to  $M$. {\it The reduction number} of $S_+$ with respect to $M$ is defined by  $$r(S_+;M) =\min\{r_{J}(S_+;M)|J\;\text{is}\;\text{a}\text{ minimal}\text{ reduction} \text{ of}\; S_+\text{ with} \text{ respect} \text{ to}\; M\}.$$

The following lemma will determine the relationship between filter-regular sequences and minimal reductions in graded rings and modules.\\\\
{\bf Lemma 2.8} (see [Tr1]){\bf.} {\it Let $S=\bigoplus_{n\ge 0}S_n$  be a finitely generated standard graded algebra  over an Artinian local ring $A$; $M=\bigoplus_{n\ge 0}M_n$   a finitely generated graded $S$-module such that $S_1$ is not contained in $\sqrt{\mathrm{Ann}_S M}$. Set $\ell=\dim M$. Assume that $J$ is a minimal reduction of $S_+$ with respect to $M.$ Then $J$ is  generated by an $S_{+}$-filter-regular sequence with respect to $M$ consisting of $\ell$ homogeneous elements of  degree 1.}\\

Denote by $r(S_+;M)$ the reduction number of $S_+$ with respect to $M$. The following  lemma makes up an important role  in establishing the relationship between mixed multiplicities of multigraded modules and the length of modules.\\\\
{\bf Lemma 2.9.} {\it Let $S=\bigoplus_{n\ge 0}S_n$  be a finitely generated standard graded algebra  over an Artinian local ring $A$; $M=\bigoplus_{n\ge 0}M_n$   a finitely generated graded $S$-module such that $S_1$ is not contained in $\sqrt{\mathrm{Ann}_S M}$. Assume that $\dim M=1$. Set $r=r(S_+;M)$ and  $M^*=\dfrac{M}{0_M:S_+^{\infty}}.$ Then 
\begin{itemize}
\item[$\mathrm{(i)}$] $l_A(M^*_n)=l_A(M^*_r)$  for all $n\ge r$.
\item[$\mathrm{(ii)}$]  $l_A(M_n)=l_A(M^*_r)$  for all large $n$.
\end{itemize}}
\begin{proof}\:  Since  $$l_A(M^*_n)=l_A\biggl[\dfrac{M_n}{(0_M:S_+^{\infty})_n}\biggl]
=l_A(M_n)-l_A[(0_M:S_+^{\infty})_n]$$ and $(0_M:S_+^{\infty})_n=0$ for all large $n$ by Remark 2.3, $l_A(M^*_n)=l_A(M_n)$ for all large $n$. 

The proof of (i):\: Since $\dim M=1,$ by Lemma 2.8,   there exists an $S_{+}$-filter-regular element $x\in S_1$ with respect to $M$ such that $xS$ is a minimal reduction of $S_+$ with respect to $M$ and $r=r_{xS}(S_+;M)$. It is clear that $M^*_n=x^{n-r}M^*_r$  for all $n\ge r$. Therefore,  $$l_A(M^*_n)=l_A(x^{n-r}M^*_r)$$ for all $n\ge r$. Since $x$ is an $S_{+}$-filter-regular element, $x$ is  non-zero-divisor in $M^*.$ It follows that  $x^{n-r}M^*_r\cong M^*_r$. Hence    $$l_A(M^*_n)=l_A(x^{n-r}M^*_r)=l_A(M^*_r)$$ for all $n\ge r$. 

The proof of (ii):\: Since $l_A(M^*_n)=l_A(M_n)$ for all large $n$ and by (i),  we immediately obtain $l_A(M_n)=l_A(M^*_r)$ for all large $n$. $\blacksquare$\\\\
\end{proof} 

\centerline{\Large\bf3.  Mixed Multiplicities of Multigraded Modules}$\\$
\hspace*{0.5cm}Basing on the properties of filter-regular sequences (Section 2), this section answers to  the question when mixed multiplicities of multigraded modules are positive and characterizes  them in terms of  lengths of modules.

 Recall that a polynomial $P(n_1,\ldots,n_d)$ is called the {\it Hilbert  polynomial of } the function  $l_A[M_{(n_1,\ldots,n_d)}]$ if 
$P(n_1,\ldots,n_d)= l_A[M_{(n_1,\ldots,n_d)}]$ for all large $n_1,\ldots,n_d.$\\
 
\hspace*{-0.6cm}{\bf Remark 3.1.} Set $\ell=\dim M^\triangle.$ Denote by $P(n_1,\ldots,n_d)$ the Hilbert polynomial of $$l_A[M_{(n_1,\ldots,n_d)}].$$ Assume that $\ell > 0.$ By [HHRT, Theorem 4.1], $\deg P(n_1,\ldots,n_d)=\dim \text{Supp}_{++}M$  and all coefficients of monomials of highest degree in $P(n_1,\ldots,n_d)$ are non-negative integers not all zero (see [Ba]). So  $\deg P(n_1,\ldots,n_d)=\deg P(n,\ldots,n).$ Since $$P(n,\ldots,n)=l_A[M_{(n,\ldots,n)}]=l_A(M^\triangle_n)$$ for all large $n$, $$\deg P(n,\ldots,n)= \dim M^\triangle-1=\ell-1.$$  Hence  $\deg P(n_1,\ldots,n_d)=\dim\text{Supp}_{++}M =\ell-1.$\\  

\hspace*{-0.6cm}{\bf Remark 3.2.} $\dim M^\triangle = \ell > 0$ if and only if $S_{(1,1,\ldots,1)}\nsubseteq \sqrt{\mathrm{Ann}_S M}.$ Indeed, $\dim M^\triangle = 0$ is equivalent to $[M^\triangle]_n = 0$ for all large $n.$ 
The latter means the same as $M_{(n,\ldots,n)} = 0$ for all large $n,$ but it is easily seen to be equivalent to $M_{(n_1,\ldots,n_d)}= 0$ for all large $n_1,\ldots,n_d.$ Hence $\dim M^\triangle = 0$ if and only if $S_{(1,1,\ldots,1)}\subseteq \sqrt{\mathrm{Ann}_S M}$ by Proposition 2.7.   

 Denote by  $L(M)$ the set of the lengths of maximal $S_{++}$-filter-regular sequences in $\bigcup _{j=1}^dS_j$ with respect to $M$. Then we have the following proposition.\\\\
{\bf Proposition 3.3.} {\it   Let $S$  be a finitely generated standard $d$-graded algebra over  an Artinian local ring $A$ and $M$  a finitely generated $d$-graded $S$-module such that $S_{(1,1,\ldots,1)}$ is not contained in $ \sqrt{\mathrm{Ann}_S M}$. Set $\ell=\dim M^\triangle.$ Assume that $$e(M;k_1,\ldots,k_d)\ne 0,$$ where $k_1,\ldots,k_d$ are non-negative integers such that $k_1+\cdots+k_d=\ell-1$.  Then the following statements hold.
\begin{itemize}
\item[$\mathrm{(i)}$] If  $k_i>0$ and $x\in S_i$ is an $S_{++}$-filter-regular element  with respect to $M$ then $$e(M;k_1,\ldots,k_d)=e(M/xM;k_1,\ldots,k_i-1,\ldots,k_d)$$ and $\dim (M/xM)^\triangle=\ell-1.$ 
\item[$\mathrm{(ii)}$] There exists an $S_{++}$-filter-regular sequence in $\bigcup _{j=1}^dS_j$ with respect to $M$ consisting of $k_1$ elements of $S_1,\ldots,k_d$ elements of $S_d$.
\item[$\mathrm{(iii)}$] $\max L(M)=\ell$.
\end{itemize}}
\begin{proof}\: Denote by $P(n_1,\ldots,n_d)$ the polynomial of $l_A[M_{(n_1,\ldots,n_d)}]$. Since $$S_{(1,1,\ldots,1)}\nsubseteq \sqrt{\mathrm{Ann}_S M},$$ by Remark 3.1 and Remark 3.2 we have $\deg P=\ell-1\ge 0.$

The proof of (i): By Remark 2.5, $$l_A[(M/xM)_{(n_1,\ldots,n_d)}]=l_A[M_{(n_1,\ldots,n_{d})}]-l_A[M_{(n_1,\ldots,n_i-1,\ldots,n_{d})}]$$ for all large $n_1,\ldots,n_d$. Denote by 
$Q(n_1,\ldots,n_d)$ the polynomial of $$l_A[(M/xM)_{(n_1,\ldots,n_d)}].$$ From the above fact, we have $$Q(n_1,\ldots,n_d)=P(n_1,\ldots,n_i,\ldots,n_d)-P(n_1,\ldots,n_i-1,\ldots,n_d).$$ Since $e(M;k_1,\ldots,k_d)\ne 0$ and $k_i>0$, we get $\deg Q =\deg P -1$ and $$e(M;k_1,\ldots,k_d)=e(M/xM;k_1,\ldots,k_i-1,\ldots,k_d).$$ By Remark 3.1,  $\deg Q=\dim (M/xM)^\triangle-1.$ Therefore $$\dim (M/xM)^\triangle=\deg Q+1.$$ Since $\deg Q =\deg P -1$ and $\deg P =\ell-1$, $\deg Q=\ell-2.$  Hence  $$ \dim (M/xM)^\triangle=\deg Q+1=(\ell-2)+1=\ell-1.$$ 

The proof of (ii): Since $S_{(1,1,\ldots,1)}\nsubseteq \sqrt{\mathrm{Ann}_S M} $, $\ell\ge 1$ by Remark 3.2. The proof is by induction on $\ell$. For $\ell=1$, the result is trivial. Assume that $\ell>1$. Since $k_1+\cdots+k_d=\ell-1>0$, there exists $k_j>0$. By Proposition 2.6, there exists an $S_{++}$-filter-regular element  $x_1\in S_j$  with respect to $M$. By (i), $$e(M;k_1,\ldots,k_d)=e(M/x_1M;k_1,\ldots,k_j-1,\ldots,k_d).$$ Since $k_1+\cdots+(k_j-1)+\cdots+k_d=\ell-2$, by the inductive assumption, there exists an $S_{++}$-filter-regular sequence $x_2,\ldots,x_{\ell-1}\in\bigcup _{j=1}^dS_j$ with respect to $M/x_1M$ consisting of $k_1$ elements of $S_1,\ldots,(k_j-1)$ elements of $S_j$, $\ldots,k_d$ elements of $S_d$. Hence $$x_1,\ldots,x_{\ell-1}\in\bigcup _{j=1}^dS_j$$ is an 
$S_{++}$-filter-regular sequence with respect to $M$ consisting of $k_1$ elements of $S_1,\ldots,k_j$ elements of $S_j,\ldots,k_d$ elements of $S_d.$

The proof of (iii): We first show that $\max L(M)\ge \ell.$ By (ii), there exists an $S_{++}$-filter-regular sequence $x_1,\ldots,x_{\ell-1}\in\bigcup _{j=1}^dS_j$ with respect to $M$ consisting of $k_1$ elements of $S_1,\ldots,k_d$ elements of $S_d$. By (i), $$e(M/(x_1,\ldots,x_{\ell-1})M;0,\ldots,0)=e(M;k_1,\ldots,k_d)\ne 0.$$ Since $e(M/(x_1,\ldots,x_{\ell-1})M;0,\ldots,0)\ne 0,$ it follows that  $$\dim [M/(x_1,\ldots,x_{\ell-1})M]^\triangle > 0.$$ Hence by Remark 3.2,    
$$S_{(1,1,\ldots,1)}\nsubseteq\sqrt{\mathrm{Ann}_S[M/(x_1,\ldots,x_{\ell-1})M]}.$$  By Proposition 2.6, there exists an $S_{++}$-filter-regular element  $x_{\ell}\in \bigcup _{j=1}^dS_j$ with respect to $M/(x_1,\ldots,x_{\ell-1})M.$ Hence $x_1,\ldots,x_{\ell}$ is an $S_{++}$-filter-regular sequence in $\bigcup _{j=1}^dS_j$ with respect to $M$. Since the length of this sequence is $\ell$,  $\max L(M)\ge \ell.$ Now, we need to show that $\max L(M)\le \ell.$ Note that $\ell=\deg P+1$. We will  prove that $\max L(M)\le \deg P+1$ by induction on $\deg P$. For each $j=1,\ldots,d$, suppose that $x_1\in S_j$ is an $S_{++}$-filter-regular element with respect to $M$. Denote by $Q_1(n_1,\ldots,n_d)$ the polynomial of $l_A[(M/x_1M)_{(n_1,\ldots,n_d)}]$. By Remark 2.5, $$l_A[(M/x_1M)_{(n_1,\ldots,n_d)}]=l_A[M_{(n_1,\ldots,n_j,\ldots,n_d)}]-l_A[M_{(n_1,\ldots,n_j-1,\ldots,n_d)}]$$ for all large $n_1,\ldots,n_d.$ Hence
$$Q_1(n_1,\ldots,n_d)=P(n_1,\ldots,n_j,\ldots,n_d)-P(n_1,\ldots,n_j-1,\ldots,n_d).$$  This implies that $\deg Q_1\le \deg P-1.$ By the inductive assumption, the length of any maximal $S_{++}$-filter-regular sequence in $\bigcup _{j=1}^dS_j$ with respect to $M/x_1M$ is not greater than $\deg Q_1 +1$. So the length of any maximal $S_{++}$-filter-regular sequence in $\bigcup _{j=1}^dS_j$ with respect to $M$ is not greater than $\deg Q_1 +2.$ Hence $$\max L(M)\le \deg Q_1 +2\le (\deg P-1)+2=\deg P+1=\ell.$$  By the above facts, we get $\max L(M)=\ell$. $\blacksquare$
\end{proof}
The following theorem will give the criteria for the positivity of  mixed multiplicities and characterize them in terms of lengths of modules.\\\\
{\bf Theorem 3.4.} {\it Let $S$ be a finitely generated standard $d$-graded algebra over  an Artinian local ring $A$ and $M$  a finitely generated $d$-graded $S$-module such that $S_{(1,1,\ldots,1)}$ is not contained in $ \sqrt{\mathrm{Ann}_S M}$. Set $\ell=\dim M^\triangle.$ Then the following statements hold.
\begin{itemize}
\item[$\mathrm{(i)}$] $e(M;k_1,\ldots,k_d)\ne 0 $\; if and only if  there exists an $S_{++}$-filter-regular sequence $x_1,\ldots,x_{\ell-1}$  with respect to  $M$ consisting of $k_1$ elements of $S_1,\ldots,k_d$ elements of  $S_d$ and 
$\dim \ovl M^{\triangle} =1$, where $\ovl M=\dfrac{M}{(x_1,\ldots,x_{\ell-1})M}$.
\item[$\mathrm{(ii)}$] Suppose that $e(M;k_1,\ldots,k_d)\ne 0 $ and $x_1,\ldots,x_{\ell-1} $ is  an $S_{++}$-filter-regular sequence  with respect to  $M$ consisting of $k_1$ elements of $S_1,\ldots,k_d$ elements of  $S_d$. Set $\ovl M=\dfrac{M}{(x_1,\ldots,x_{\ell-1})M}\; and\; r=r(S^{\triangle}_+;\ovl M^{\triangle})$ the reduction number of $S^{\triangle}_+$ with respect to  $\ovl M^{\triangle}.$  Then $ e(M;k_1,\ldots,k_d)=l_A\biggl[\dfrac{\ovl M^\triangle_n}{(0_{\ovl M^\triangle}:S_+^{\triangle\infty})_n}\biggl]\;\text{for\; all\;} n\ge r.$
\end{itemize}}

\begin{proof}\: The proof of (i):  First, we prove the necessary condition.  By Proposition 3.3(ii), there exists an $S_{++}$-filter-regular sequence   $x_1,\ldots,x_{\ell-1}$ with respect to $M$ consisting of $k_1$ elements of $S_1,\ldots,k_d$ elements of  $S_d$. Applying Proposition 3.3(i), by induction on $\ell$ we have $\dim \ovl M^{\triangle}=\ell-(\ell-1)=1.$ Hence we obtain the necessary condition. 
 
We turn to the proof of sufficient condition. The proof is by induction on $\ell\ge 1$. Denote by $P(n_1,\ldots,n_d)$ the polynomial of $l_A[M_{(n_1,\ldots,n_d)}]$. Note that by Remark 3.1, $\deg P=\ell-1$. For $\ell=1$, we have $\deg P=0$. So $P(n_1,\ldots,n_d)$ is a constant not zero. Hence $$e(M;k_1,\ldots,k_d)=e(M;0,\ldots,0)=P(n_1,\ldots,n_d)\ne 0.$$ If $\ell>1,$ since $k_1+\cdots+k_d=\ell-1>0$, there exists $1\le i\le d$ such that $k_i>0.$ And without loss of general, we may assume that $x_1\in S_i$. 
  Set $N = M/x_1M.$ Denote by $Q(n_1,\ldots,n_d)$ the polynomial of $l_A[N_{(n_1,\ldots,n_d)}]$. By the proof of Proposition 3.3(i), we have $$Q(n_1,\ldots,n_d)=P(n_1,\ldots,n_i,\ldots,n_d)-P(n_1,\ldots,n_i-1,\ldots,n_d).$$ This implies that $\deg Q\le \deg P-1=\ell-2.$ Since $\dim \ovl M^{\triangle}=1$, it follows that  $S_{(1,1,\ldots,1)} \nsubseteq\sqrt{\mathrm{Ann}_S \ovl M }$ by Remark 3.2. By Proposition 2.6, there exists an $S_{++}$-filter-regular element $x_\ell$ in $\bigcup_{j=1}^dS_j$ with respect to  $\ovl M$. By Proposition 3.3(iii), $x_1,\ldots,x_\ell$ is a maximal $S_{++}$-filter-regular sequence in $\bigcup_{j=1}^dS_j$ with respect to $M$. 
It is clear that $x_2,\ldots,x_\ell$ is a maximal $S_{++}$-filter-regular sequence  in $\bigcup_{j=1}^dS_j$ with respect to $N.$  Denote by $L(N)$ the set of the lengths of maximal $S_{++}$-filter-regular sequences in $\bigcup _{j=1}^dS_j$ with respect to $N.$
We immediately get $ \max L(N)\ge \ell-1$. By Remark 3.1 and Proposition 3.3(iii), 
$ \max L(N)=\dim N^{\triangle}=\deg Q+1.$ Hence $\deg Q\ge \ell-2.$ By the above facts, we get  $\deg Q=\ell-2=\deg P-1.$ 
From this equality and $$Q(n_1,\ldots,n_d)=P(n_1,\ldots,n_i,\ldots,n_d)-P(n_1,\ldots,n_i-1,\ldots,n_d),$$ we have $$e(M;k_1,\ldots,k_d)=e(N;k_1,\ldots,k_i-1,\ldots,k_d).$$
 It is a plain fact  that $x_2,\ldots,x_{\ell-1}$ is an $S_{++}$-filter-regular sequence of the length  $\ell-2$ with respect to $N$ consisting of $k_1$ elements of $S_1,\ldots,k_i-1$ elements of  $S_i,\ldots,k_d$ elements of  $S_d.$  Set $\ovl N=N/(x_2,\ldots,x_{\ell-1})N.$  It can be verified that $\ovl M\cong\ovl N.$ Hence 
 $$\dim \ovl N^{\triangle}=\dim \ovl M^{\triangle} =1.$$ Applying the inductive assumption for $N,$ we have $e(N;k_1,\ldots,k_i-1,\ldots,k_d) \ne 0.$ So $e(M;k_1,\ldots,k_d)\ne 0$.

The proof of (ii):  Applying  Proposition 3.3(i),  by induction on $\ell$  we get $$ 0 \ne e(M;k_1,\ldots,k_d)= e(\ovl M;0,\ldots,0)$$ and $\dim \ovl M^{\triangle}=\ell-(\ell-1)=1.$  By Remark 3.1, $l_A[\ovl M_{(n_1,\ldots,n_d)}]$ is a polynomial of degree $0$ for all large $n_1,\ldots,n_d.$ Hence $e(\ovl M;0,\ldots,0) = l_A[\ovl M_{(n_1,\ldots,n_d)}]$ for all large $n_1,\ldots,n_d.$ By taking $n_1=\cdots=n_d=n$, where $n$ is a sufficiently large positive integer, we have $e(\ovl M;0,\ldots,0)=l_A[\ovl M_{(n,\ldots,n)}]=l_A[\ovl M^{\triangle}_n]$ for all large $n$.  Consequently,  $$e(M;k_1,\ldots,k_d)=e(\ovl M;0,\ldots,0)=l_A[\ovl M^{\triangle}_n]$$ for all large $n$.  By Lemma 2.9,  $l_A[\ovl M^{\triangle}_n]=l_A\biggl[\dfrac{\ovl M^\triangle_r}{(0_{\ovl M^\triangle}:S_+^{\triangle \infty})_r}\biggl]$ for all large $n$ and $$l_A\biggl[\dfrac{\ovl M^\triangle_n}{(0_{\ovl M^\triangle}:S_+^{\triangle \infty})_n}\biggl]=l_A\biggl[\dfrac{\ovl M^\triangle_r}{(0_{\ovl M^\triangle}:S_+^{\triangle \infty})_r}\biggl]$$  for all $n\ge r.$   Hence $ e(M;k_1,\ldots,k_d)=l_A\biggl[\dfrac{\ovl M^\triangle_n}{(0_{\ovl M^\triangle}:S_+^{\triangle \infty})_n}\biggl]$ for all $n\ge r.$ $\blacksquare$$\\\\$ 
\end{proof} 

  From the proofs of Proposition 3.3 and Theorem 3.4. we immediately give the follwing result.\\\\    
{\bf Corollary 3.5.} {\it Let $S$ be a finitely generated standard $d$-graded algebra over  an Artinian local ring $A$ and $M$  a finitely generated $d$-graded $S$-module such that $S_{(1,1,\ldots,1)}$ is not contained in $ \sqrt{\mathrm{Ann}_S M}$. Set $\ell=\dim M^\triangle.$ Then the following statements hold.
\begin{itemize}
\item[$\mathrm{(i)}$] If $e(M;k_1,\ldots,k_d)\ne 0 $\; then for any  $S_{++}$-filter-regular sequence $x_1,\ldots,x_{n}$  with respect to  $M$ consisting of $m_1 \le k_1$ elements of $S_1,\ldots,m_d \le k_d$ elements of  $S_d$ we have $\dim \ovl M^{\triangle} = \ell -n$ and $e(\ovl M;k_1-m_1,\ldots,k_d-m_d)\ne 0 $, where $\ovl M=\dfrac{M}{(x_1,\ldots,x_{n})M}.$ 
\item[$\mathrm{(ii)}$] If there exists an $S_{++}$-filter-regular sequence  $x_1,\ldots,x_{n}$ with respect to  $M$ consisting of $m_1 \le k_1$ elements of $S_1,\ldots,m_d \le k_d$ elements of  $S_d$ such that $\dim \ovl M^{\triangle} = \ell -n$ and $e(\ovl M;k_1-m_1,\ldots,k_d-m_d)\ne 0 $, where $\ovl M=\dfrac{M}{(x_1,\ldots,x_{n})M},$ then $e(M;k_1,\ldots,k_d)\ne 0 .$
\item[$\mathrm{(iii)}$] Suppose that $e(M;k_1,\ldots,k_d)\ne 0 $ and $x_1,\ldots,x_{n} $ is  an $S_{++}$-filter-regular sequence  with respect to  $M$ consisting of $m_1\le k_1$ elements of $S_1,\ldots,m_d \le k_d$ elements of   $S_d$.  Then $ e(M;k_1,\ldots,k_d)= e(\ovl M;k_1-m_1,\ldots,k_d-m_d)$ and $\dim \ovl M^{\triangle} = \ell -n,$  where $\ovl M=\dfrac{M}{(x_1,\ldots,x_{n})M}.$
  
\end{itemize}}

\noindent
{\bf Example 3.6:} Let  $(R, \frak n)$  be  a  Noetherian   local ring with  maximal ideal $\frak{n},$  and  an ideal  $\frak n$-primary  $J.$  Set $B = R/J,$ $B$ is an Artinian local ring.
 Let $$S = B\big[X_1, X_2,\ldots,X_t\big]$$ be the  ring of polynomials in  $t$ indeterminates  
$X_1, X_2,\ldots,X_t$  with coefficients in $B.$ Then $\dim S = t$ and $S$ is a finitely generated standard  graded algebra over $B.$ Since $X_1, X_2,\ldots,X_t$ is a regular sequence, $X_1, X_2,\ldots,X_t$ is   
an $S_+$-filter-regular sequence.  It is a plain fact that the reduction number of $S_+$ is $r = r(S_+) = 0$ and $\dim S/(X_1, X_2,\ldots,X_{t-1}) = 1,$ and 
$S/(X_1, X_2,\ldots,X_{t-1})\cong B[X_t].$   
Hence by Theorem 3.4, $e(S;t-1) \ne 0$ and $e(S;t-1)= l_B(B)= l_R(R/J).$ \\ 

\noindent
{\bf Example 3.7:} Let $k$ be an infinite field and  let $x_1,x_2,x_3,y_1,y_2,y_3,z_1,z_2,z_3$ be indeterminates. Let $R=k[x_1,x_2,x_3,y_1,y_2,y_3,z_1,z_2,z_3]$ be a finitely generated standard 3-graded algebra over $k$ with   $\deg x_i=(1,0,0),\deg y_i=(0,1,0), \deg z_i=(0,0,1), $$i=1,2,3.$  Set  $I=(x_1,y_1,z_1)\cap (x_1,x_2)\cap (y_1,y_2)\cap (z_1,z_2)$ and  $S=R/I.$ Then  $S$ is a finitely generated standard 3-graded algebra over $k$ and  $\dim S=7.$  Denote by  $f(n_1,n_2,n_3)$ the polynomial such that   $f(n_1,n_2,n_3)=l_k[S_{(n_1,n_2,n_3)}]$ for all large $n_1,n_2,n_3.$ We have   $$\dim S/S_{(1+)}=\dim S/S_{(2+)}=\dim S/S_{(3+)}=6=\dim S-1$$  and  $\dim \dfrac{S}{S_{(i+)}+S_{(j+)}}=3<\dim S-2$  for all $i,j=1,2,3,i\ne j.$ Hence by [HHRT, Theorem 4.3], $\deg f(n_1,n_2,n_3)=\dim S-3=7-3=4$ and $\ell=\dim S^\triangle=5.$ For $x\in R,$ denote by $\bar x$ the image of $x$ in $S.$
Direct coputation shows that $\bar x_3,\bar x_2,\bar y_3,\bar y_2$ is an $S_{++}$-filter-regular sequence of  $S$ consisting of  2 elements of  $S_1$ and  2 elements of  $S_2$.  It can be verified that   $$\begin{array}{ll}\dfrac{S}{(\bar x_3):S_{++}^\infty}&\cong \dfrac{R}{[(x_3)+I]:R_{++}^\infty}= \dfrac{R}{(x_3)+(x_1,y_1,z_1)\cap  (y_1,y_2)\cap (z_1,z_2)}

\vspace{6pt}\\&= \dfrac{R}{(x_1,x_3,y_1,z_1)\cap  (x_3,y_1,y_2)\cap (x_3,z_1,z_2)},

\vspace{6pt}\\\dfrac{S}{(\bar x_3,\bar x_2):S_{++}^\infty}&\cong \dfrac{R}{[(x_2,x_3)+I]:R_{++}^\infty}= \dfrac{R}{(x_2,x_3)+  (y_1,y_2)\cap (z_1,z_2)}

\vspace{6pt}\\&= \dfrac{R}{  (x_2,x_3,y_1,y_2)\cap (x_2,x_3,z_1,z_2)},

\vspace{6pt}\\\dfrac{S}{(\bar x_3,\bar x_2,\bar y_3):S_{++}^\infty}&\cong \dfrac{R}{[(x_2,x_3,y_3)+I]:R_{++}^\infty}= \dfrac{R}{(x_2,x_3,y_3,z_1,z_2)}

\vspace{6pt}\\&\cong k[x_1,y_1,y_2,z_3]\end{array}$$ and  $\dfrac{S}{(\bar x_3,\bar x_2,\bar y_3,\bar y_2):S_{++}^\infty}\cong k[x_1,y_1,z_3].$ Since  $$\dim \biggl[\dfrac{S}{(\bar x_3,\bar x_2,\bar y_3,\bar y_2)}\biggl]^\triangle=\dim \biggl[\dfrac{S}{(\bar x_3,\bar x_2,\bar y_3,\bar y_2):S_{++}^\infty}\biggl]^\triangle=\dim k[x_1,y_1,z_3]^\triangle=1,$$ by Theorem 3.4 we get $$\begin{array}{ll}0\ne e(S;2,2,0)&=e(\dfrac{S}{(\bar x_3,\bar x_2,\bar y_3,\bar y_2)};0,0,0)=e(\dfrac{S}{(\bar x_3,\bar x_2,\bar y_3,\bar y_2):S_{++}^\infty};0,0,0)\\\\&=e(k[x_1,y_1,z_1];0,0,0)=1.\end{array}$$ By symmetry, we also have  $e(S;2,0,2)=e(S;0,2,2)=1.$
It can be verified that  $\bar x_3,\bar x_2,\bar x_1$ is also  an $S_{++}$-filter-regular sequence consisting of  3 elements of  $S_1.$  Since  $\dfrac{S}{(\bar x_3,\bar x_2,\bar x_1):S_{++}^\infty}=0,$ $\bar x_3,\bar x_2,\bar x_1$ is a maximal  $S_{++}$-filter-regular sequence. Hence $e(S;3,1,0)= e(S;4,0,0)=0.$ Indeed, if $e(S;3,1,0)\ne 0$ or $e(S;4,0,0)\ne 0$ then by Corollary 3.5, $\bar x_3,\bar x_2,\bar x_1$ is not a maximal $S_{++}$-filter-regular sequence. By symmetry, we get $$\begin{array}{ll}&e(S;3,1,0)=e(S;1,3,0)=e(S;3,0,1)=e(S;1,0,3)=e(S;0,3,1)\\\\&=e(S;0,1,3)=e(S;4,0,0)=e(S;0,4,0)=e(S;0,0,4)=0.\end{array}$$
Upon simple computation, we show that    $\bar x_3,\bar x_2,\bar y_3,\bar z_3$ is an  $S_{++}$-filter-regular sequence of $S$ consisting of  2 elements of $S_1$, 1 element of  $S_2$ and  1 element of  $S_3$. Since $\dfrac{S}{(\bar x_3,\bar x_2,\bar y_3,\bar z_3):S_{++}^\infty}=0,$ $e(S;2,1,1)=0.$ Indeed, if $e(S;2,1,1)\ne 0$ then by Theorem 3.4(ii), $$0 \ne e(S;2,1,1)= l_k\biggl[\dfrac{S}{(\bar x_3,\bar x_2,\bar y_3,\bar z_3):S_{++}^\infty}\biggl]_{(n,\ldots,n)}$$ for all large $n.$ 
By symmetry, we get $$e(S;2,1,1)=e(S;1,2,1)=e(S;1,1,2)=0.$$ From the above  facts, we also obtain
 $$e(S)=\sum_{k_1+k_2+k_3=4}e(S;k_1,k_2,k_3)=3$$  by [HHRT, Theorem 4.3].\\\\

\centerline{\Large\bf4. Applications}$\\\\$
\hspace*{0.5cm} This section is devoted to the discussion of applications of Section 3. We deal with the positivity and the relationship between  mixed multiplicities of ideals  and Hilbert-Samuel multiplicities.   

 \vskip 0.2cm

 Let  $(R, \frak n)$  be  a  Noetherian   local ring of Krull dimension $\dim R = d >0$ with maximal ideal $\frak{n}$ and infinite residue field $k = R/\frak{n},$  and  an   $\frak n$-primary ideal  $J$, and $I_1,\ldots, I_s$ ideals of $R$ such that $I= I_1\cdots I_s$ is non-nilpotent.  
Remember that 

\vspace{3pt}
$$  T  =\bigoplus_{n, n_1,\ldots,n_s\ge 0}\dfrac{J^nI_1^{n_1}\cdots I_s^{n_s}}{J^{n+1}I_1^{n_1}\cdots I_s^{n_s}}$$
\vspace{6pt}is a finitely generated standard multigraded algebra over  Artinian local ring $R/J,$ and mixed multiplicities of $T$ are 
mixed multiplicities of ideals $J,I_1,\ldots, I_s$ (see [Ve] or [HHRT]). 
By [Vi], $\dim T^\triangle = \dim R/0: I^\infty.$ Set $$\dim R/0: I^\infty = q;\; e(T;k_0, k_1,\ldots,k_s)  = e(J^{[k_0+1]},I_1^{[k_1]},\ldots,I_s^{[k_s]}; R)$$ $(k_0+k_1+\cdots+k_s = q-1).$ Then $e(J^{[k_0+1]},I_1^{[k_1]},\ldots,I_s^{[k_s]}; R)$ is called the   mixed multiplicity of ideals $J,I_1,\ldots,I_s$ 
of type $(k_0,k_1,\ldots,k_s).$  

 \vskip 0.2cm
Using different sequences, one transmuted  mixed multiplicities into  Hilbert-Samuel multiplicities, for instance: Risler and Teissier in 1973 [Te] by superficial sequences; Rees in 1984 [Re] by joint reductions; Viet in 2000 [Vi] by (FC)-sequences; Trung and Verma in 2007 [TV] by $(\varepsilon_1,\ldots,\varepsilon_m)$-superficial sequences. 
\vskip 0.2cm
  
Our approach is based on the results in Section 3. Assume that $ x \in I_i$ and $x^*$ is the image of $x$ in $ T.$  We need to choose $x$  satisfying the following properties: 
\begin{itemize}
 
 \vspace{3pt}
 
 \item[\rm (i):]$x^*$  is an $T_{++}$-filter-regular element of $ T.$

\vspace{3pt}
\item [\rm (ii):]$[T/x^*T]_{(m, m_1,\ldots,m_s)}\cong \bigg[\bigoplus_{n, n_1,\ldots,n_s\ge 0}\dfrac{J^nI_1^{n_1}\cdots I_s^{n_s}(R/xR)}{J^{n+1}I_1^{n_1}\cdots I_s^{n_s}(R/xR)}\bigg]_{(m, m_1,\ldots,m_s)}$ for all $m, m_1,\ldots,m_s\gg0.$  
\end{itemize}

Under our point of view in this paper, both conditions are necessary to express mixed multiplicities of ideals in terms of Hilbert-Samuel multiplicity. Now we will characterize these conditions by equations of ideals.

\vskip 0.2cm

\hspace*{-0.6cm}{\bf Remark 4.1.} Let $x \in I_i$ and $x^*$ the image of $x$ in $T.$ 
It can be verified that $x^*$  is an $T_{++}$-filter-regular element of $ T$ if and only if      
 $$(J^{n+1}I_1^{n_1}\cdots I_i^{n_i+1}\cdots I_s^{n_s}:x) \bigcap J^{n}I_1^{n_1}\cdots I_s^{n_s} = J^{n+1}I_1^{n_1}\cdots I_s^{n_s}$$ for all $n, n_1,\ldots,n_s\gg0,$ and 
$$[T/x^*T]_{(m, m_1,\ldots,m_s)}\cong\dfrac{J^mI_1^{m_1}\cdots I_s^{m_s}}{J^{m+1}I_1^{m_1}\cdots I_s^{m_s}+xJ^{m}I_1^{m_1}\cdots I_i^{m_i-1}\cdots I_s^{m_s}},$$
and  $$\begin{array}{l}\bigg[\bigoplus_{n, n_1,\ldots,n_s\ge 0}\dfrac{J^nI_1^{n_1}\cdots I_s^{n_s}(R/xR)}{J^{n+1}I_1^{n_1}\cdots I_s^{n_s}(R/xR)}\bigg]_{(m, m_1,\ldots,m_s)}\cong \dfrac{J^mI_1^{m_1}\cdots I_s^{m_s}+(x)}{J^{m+1}I_1^{m_1}\cdots I_s^{m_s}+(x)}\\

\vspace{6pt}\cong \dfrac{J^mI_1^{m_1}\cdots I_s^{m_s}}{J^{m+1}I_1^{m_1}\cdots I_s^{m_s}+(x)\bigcap J^mI_1^{m_1}\cdots I_s^{m_s}}.\end{array}$$
Hence  $$[T/x^*T]_{(m, m_1,\ldots,m_s)} \cong \bigg[\bigoplus_{n, n_1,\ldots,n_s\ge 0}\dfrac{J^nI_1^{n_1}\cdots I_s^{n_s}(R/xR)}{J^{n+1}I_1^{n_1}\cdots I_s^{n_s}(R/xR)}\bigg]_{(m, m_1,\ldots,m_s)}$$ for all $m, m_1,\ldots,m_s\gg0$
if and only if $$J^{m+1}I_1^{m_1}\cdots I_s^{m_s}+(x)\bigcap J^mI_1^{m_1}\cdots I_s^{m_s}=J^{m+1}I_1^{m_1}\cdots I_s^{m_s}+xJ^{m}I_1^{m_1}\cdots I_i^{m_i-1}\cdots I_s^{m_s}$$ for all $m, m_1,\ldots,m_s\gg0$
or $$(x)\bigcap J^mI_1^{m_1}\cdots I_s^{m_s}\equiv xJ^{m}I_1^{m_1}\cdots I_i^{m_i-1}\cdots I_s^{m_s}\;(\mathrm{mod }\:J^{m+1}I_1^{m_1}\cdots I_s^{m_s})$$ for all $m, m_1,\ldots,m_s\gg0.$ Note that if
$$[T/x^*T]_{(m, m_1,\ldots,m_s)} \cong \bigg[\bigoplus_{n, n_1,\ldots,n_s\ge 0}\dfrac{J^nI_1^{n_1}\cdots I_s^{n_s}(R/xR)}{J^{n+1}I_1^{n_1}\cdots I_s^{n_s}(R/xR)}\bigg]_{(m, m_1,\ldots,m_s)}$$ for all $m, m_1,\ldots,m_s\gg0$ then   
$$\dim (T/x^*T)^\triangle = \dim R/(x): I^\infty.$$ Denote by  $\bar{J}, \bar{I}_1,\ldots,\bar{I}_s$ the images of $J, I_1,\ldots,I_s$ in $R/(x).$ Set $\dim (T/x^*T)^\triangle = t.$ Then we have  
$$e(T/x^*T;k_0, k_1,\ldots,k_s)  = e(\bar{J}^{[k_0+1]},\bar{I}_1^{[k_1]},\ldots,\bar{I}_s^{[k_s]}; R/xR)$$ $(k_0+k_1+\cdots+k_s = t-1).$ 
 \vskip 0.2cm

The above comments are the reason for using the following sequences.
 \vskip 0.2cm
\noindent {\bf Definition 4.2.} {\it Let $\frak I$ be an ideal such that $\frak II_1\cdots I_s$ is non nilpotent.  An element $x \in R$ is called a superficial element of $(I_1,\ldots, I_s)$  with respect to $\frak I = I_0$ if there exists $i \in \{0, 1, \ldots, s\}$ such that $x \in I_i$ and 
 
$\rm(i):$ $({\frak I}I_0^{n_0}\cdots I_i^{n_i+1}\cdots I_s^{n_s}:x) \bigcap I_0^{n_0}\cdots I_s^{n_s} = {\frak I}I_0^{n_0}\cdots I_i^{n_i}\cdots I_s^{n_s}$ for all\\\hspace*{1.7cm} $n_0, n_1,\ldots,n_s\gg0.$

$\rm (ii):$ $(x)\bigcap {I_0}^{n_0} \cdots I_i^{n_i+1}\cdots I_s^{n_s} 
= x{I_0}^{n_0}\cdots I_i^{n_i}\cdots I_s^{n_s}$
for $n_0,n_1,\ldots,n_s\gg0.$\\  
Let $x_1, \ldots, x_t$ be a sequence in $R$. For each $i = 0, 1, \ldots, t - 1 $, set $R_i = R/(x_1, \ldots, x_{i})$, $\bar {\frak I} = {\frak I}R_i , \bar{I}_j = I_jR_i$, $\bar{x}_{i + 1}$ the image of $x_{i + 1}$ in $R_i$. Then 
$x_1, \ldots, x_t$ is called a superficial sequence  of $(I_1,\ldots, I_s)$ with respect to $\frak I$ if $\bar{x}_{i + 1}$ is a superficial element  of $(\bar{I}_1,\ldots, \bar{I}_s)$  with respect to $\bar {\frak I}$ for all $i = 0, 1, \ldots, t - 1$.}

\vspace*{6pt} 

\hspace*{-0.5cm}{\bf Remark 4.3.} 
  Assume that $x \in I_i$ and $x^*$ is the image of $x$ in $ T$. Then 
 \begin{enumerate}
\item[\rm (i)]  $x^*$  is an $T_{++}$-filter-regular element of $ T$ if and only if      
 $x$ satisfies the condition (i) of a superficial element of $(I_1,\ldots, I_s)$ with respect to $J.$ 

\item[\rm (ii)] Trung and Verma [TV] called  that   
 $x$ is an $i$-superficial element for $J, I_1,\ldots,I_s$ if $x \in I_i$ and the image $x^{**}$ of $x$ in 
$ \frak T  =\bigoplus_{n, n_1,\ldots,n_s\ge 0}\dfrac{J^nI_1^{n_1}\cdots I_s^{n_s}}{J^{n+1}I_1^{n_1+1}\cdots I_s^{n_s+1}}$
is an ${\frak T}_{++}$-filter-regular element in $\frak T$, i.e., 
$$(J^{n+1}I_1^{n_1+1}\cdots I_i^{n_i+2}\cdots I_s^{n_s+1}:x) \bigcap J^{n}I_1^{n_1}\cdots I_s^{n_s} = J^{n+1}I_1^{n_1+1}\cdots I_i^{n_i+1}\cdots I_s^{n_s+1}$$
for $n, n_1,\ldots,n_s\gg0.$ And  if  $\varepsilon_1,\ldots,\varepsilon_m$ is a non-decreasing sequence of indices with $1\le \varepsilon_i\le s,$ then a sequence $x_1,\ldots,x_m$ is an  $(\varepsilon_1,\ldots,\varepsilon_m)$-superficial sequence for $J, I_1,\ldots,I_s$ if for $i= 1,\ldots,m$, $\bar{x}_i$ is an $\varepsilon_i$-superficial element for $\bar{J}, \bar{I}_1,\ldots,\bar{I}_s$, where $\bar{x}_i,\bar{J}, \bar{I}_1,\ldots,\bar{I}_s$ are the images of $x_i,J, I_1,\ldots,I_s$ in 
$$R/(x_1,\ldots,x_{i-1}).$$
Let $x$ be an $i$-superficial element for $J, I_1,\ldots,I_s.$ 
We need to show that  
$$(J^{n+1}I_1^{n_1}\cdots I_i^{n_i+1}\cdots I_s^{n_s}:x) \bigcap J^{n}I_1^{n_1}\cdots I_s^{n_s} = J^{n+1}I_1^{n_1}\cdots I_s^{n_s}$$ for $n, n_1,\ldots,n_s\gg0.$
Indeed, since 
$$(J^{n+1}I_1^{n_1+1}\cdots I_i^{n_i+2}\cdots I_s^{n_s+1}:x) \bigcap J^{n}I_1^{n_1}\cdots I_s^{n_s} = J^{n+1}I_1^{n_1+1}\cdots I_i^{n_i+1}\cdots I_s^{n_s+1}$$ for $n, n_1,\ldots,n_s\gg0,$ 
$$(J^{n+1}I_1^{n_1}\cdots I_i^{n_i+1}\cdots I_s^{n_s}:x) \bigcap J^{n}I_1^{n_1-1}\cdots I_s^{n_s-1} = J^{n+1}I_1^{n_1}\cdots I_s^{n_s}$$ for $n, n_1,\ldots,n_s\gg0.$
Hence $$\begin{array}{ll}&(J^{n+1}I_1^{n_1}\cdots I_i^{n_i+1}\cdots I_s^{n_s}:x) \bigcap J^{n+1}I_1^{n_1}\cdots I_s^{n_s} \\&=
[(J^{n+1}I_1^{n_1}\cdots I_i^{n_i+1}\cdots I_s^{n_s}:x)\bigcap J^{n}I_1^{n_1-1}\cdots I_s^{n_s-1}] \bigcap J^{n+1}I_1^{n_1}\cdots I_s^{n_s}  \\&= J^{n+1}I_1^{n_1}\cdots I_s^{n_s}\end{array}$$ 
for $n, n_1,\ldots,n_s\gg0.$ So $x$ satisfies the condition (i) of a superficial element of $(I_1,\ldots, I_s)$ with respect to $J.$ 
 On the other hand $x$ satisfies also the condition (ii) of a superficial element of $(I_1,\ldots, I_s)$ with respect to $J$ 
 by [TV]. Consequently, $x$ is a superficial element of $(I_1,\ldots, I_s)$ with respect to $J.$   

 \item[\rm (iii)] Note that since $I=I_1\cdots I_s$ is non-nilpotent,  
  for each $i=1,\ldots,s$, there exists an $i$-superficial element for $J, I_1,\ldots,I_s$ by [TV]. Hence in the case that $I$ non nilpotent, the  existence of superficial elements of $(I_1,\ldots, I_s)$ with respect to $J$ is universal  by (ii).  
\end{enumerate}

By Remark 4.1 and Remark 4.3,  one can replace ``filter-regular elements" in Corollary 3.5 by ``superficial elements" and get a version of Corollary 3.5 as follows:\\       

\noindent {\bf Proposition 4.4.} {\it Let $(R, \frak{n})$ denote a Noetherian local ring with maximal ideal 
$\mathfrak{n}$, infinite residue $k = R/\mathfrak{n},$
  and  an $\frak n$-primary ideal  $J=I_0$, and $I_1,\ldots, I_s$   ideals of $R$ such that $I = I_1\cdots I_s$ is  non nilpotent. Then the following statements hold.
\begin{enumerate}
\item[\rm (i)] If $e(J^{[k_0+1]}, I_1^{[k_1]},\ldots, I_s^{[k_s]}; R) \not = 0$ then for any superficial sequence  
$x_1,  \ldots, x_n$ of\\ 
$(I_1,\ldots, I_s)$ with respect to $J$ 
consisting of $m_0 \le k_0$ elements of $I_0$,\dots, $m_s \le k_s$ elements of $I_s$ we have $\dim R/(x_1,  \ldots, x_n): I^\infty = \dim R/0: I^\infty -n$ and $e(\bar{J}^{[k_0+1-m_0]},\bar{I}_1^{[k_1-m_1]},\ldots, \bar{I}_s^{[k_s-m_s]}; \ovl R) \not = 0,$  where  $\bar{J}, \bar{I}_1,\ldots,\bar{I}_s$ are the images of $J, I_1,\ldots,I_s$ in $\ovl{R} = R/(x_1, \ldots, x_n).$   
\item[\rm (ii)] If there exists a superficial sequence 
$x_1,  \ldots, x_n$ of $(I_1,\ldots, I_s)$ with respect to $J$ 
consisting of $m_0 \le k_0$ elements of $I_0$,\dots, $m_s \le k_s$ elements of $I_s$ such that $\dim R/(x_1,\ldots, x_n): I^\infty = \dim R/0: I^\infty -n$ and 

\centerline{$e(\bar{J}^{[k_0+1-m_0]},\bar{I}_1^{[k_1-m_1]},\ldots, \bar{I}_s^{[k_s-m_s]}; \ovl R) \not = 0,$}  where  $\bar{J}, \bar{I}_1,\ldots,\bar{I}_s$ are the images of $J, I_1,\ldots,I_s$ in $\ovl{R} = R/(x_1, \ldots, x_n),$ then $e(J^{[k_0+1]}, I_1^{[k_1]},\ldots, I_s^{[k_s]}; R) \not = 0.$  
\item[\rm (iii)] Suppose that $e(J^{[k_0+1]}, I_1^{[k_1]},\ldots, I_s^{[k_s]}; R) \not= 0$ and $x_1,  \ldots, x_n$ 
is a superficial sequence  of $(I_1,\ldots, I_s)$ with respect to $J$ 
  consisting of $m_0 \le k_0$ elements of $I_0$,\dots, $m_s \le k_s$ elements of $I_s$. Denote by $\ovl{R} = R/(x_1, \ldots, x_n)$ and $\bar{J}, \bar{I}_1,\ldots,\bar{I}_s$ the images of $J, I_1,\ldots,I_s$ in $\ovl R.$ Then
$$e(J^{[k_0+1]}, I_1^{[k_1]},\ldots, I_s^{[k_s]}; R) = e(\bar{J}^{[k_0+1-m_0]},\bar{I}_1^{[k_1-m_1]},\ldots, \bar{I}_s^{[k_s-m_s]}; \ovl R).$$
\end{enumerate}}

Recall that if $I$ is non nilpotent,  
$e(J^{[k_0+1]},I_1^{[0]},\ldots,I_s^{[0]}; R)= e(J;R/0: I^\infty )$ by [Vi].
By Proposition 4.4, $e(J^{[k_0+1]}, I_1^{[k_1]},\ldots, I_s^{[k_s]}; R) \not = 0$ if  and only if  there exists a superficial sequence 
$x_1,  \ldots, x_n$ of $(I_1,\ldots, I_s)$ with respect to $J$ 
 consisting of $k_1$ elements of $I_1$, \dots, $k_s$ elements of $I_s$ such that $$\dim R/(x_1,  \ldots, x_n): I^\infty = \dim R/0: I^\infty -n$$ and $e(\bar{J}^{[k_0+1]},\bar{I}_1^{[0]},\ldots, \bar{I}_s^{[0]}; \ovl R) \not = 0,$  where  $\bar{J}, \bar{I}_1,\ldots,\bar{I}_s$ are the images of $J, I_1,\ldots,I_s$ in $\ovl{R} = R/(x_1, \ldots, x_n).$ But since $n < \dim R/0: I^\infty$, $\dim R/(x_1,  \ldots, x_n): I^\infty > 0.$ Hence $\bar I = \bar{I}_1\cdots\bar{I}_s$ is non nilpotent. Thus $$e(\bar{J}^{[k_0+1]},\bar{I}_1^{[0]},\ldots, \bar{I}_s^{[0]}; \ovl R) = e(J;R/(x_1,  \ldots, x_n): I^\infty) \not = 0.$$ This fact follows that  $e(J^{[k_0+1]}, I_1^{[k_1]},\ldots, I_s^{[k_s]}; R) \not = 0$ if  and only if  there exists a superficial sequence 
$x_1,\ldots, x_n$ of $(I_1,\ldots, I_s)$ with respect to $J$ 
consisting of $k_1$ elements of $I_1$, \dots, $k_s$ elements of $I_s$ such that $$\dim R/(x_1,\ldots, x_n): I^\infty = \dim R/0: I^\infty -n.$$        
Then as an immediate consequence of Proposition 4.4,  we obtain the result as follows.

 \vskip 0.2cm
\noindent {\bf Theorem 4.5.} {\it Let $(R, \frak{n})$ denote a Noetherian local ring with maximal ideal 
$\mathfrak{n}$, infinite residue $k = R/\mathfrak{n},$
  and  an  $\frak n$-primary ideal $J$, and $I_1,\ldots, I_s$   ideals of $R$ such that $I = I_1\cdots I_s$ is non nilpotent. Then the following statements hold.
\begin{enumerate}
\item[\rm (i)] $e(J^{[k_0+1]}, I_1^{[k_1]},\ldots, I_s^{[k_s]}; R) \not = 0$ if  and only if there exists a superficial sequence  
$x_1,  \ldots, x_t$ of $(I_1,\ldots, I_s)$ with respect to $J$ 
consisting of $k_1$ elements of $I_1$,\ldots, $k_s$ elements of $I_s$ and $\dim R/(x_1,  \ldots, x_t): I^\infty = \dim R/0: I^\infty -t.$ 
\item[\rm (ii)] Suppose that $e(J^{[k_0+1]}, I_1^{[k_1]},\ldots, I_s^{[k_s]}; R) \not= 0$ and $x_1,\ldots, x_t$ 
is a superficial sequence of $(I_1,\ldots, I_s)$ with respect to $J$ 
 consisting of $k_1$ elements of $I_1$,\dots, $k_s$ elements of $I_s$. Set $\ovl{R} = R/(x_1, \ldots, x_t):I^\infty$. Then
$$e(J^{[k_0+1]}, I_1^{[k_1]},\ldots, I_s^{[k_s]}; R) = e(J;\ovl{R}).$$
\end{enumerate}}

\vskip 0.2cm
\hspace*{-0.6cm}{\bf Remark 4.6.} Return to Theorem 1.4 in [TV], assume that $x_1,\ldots,x_m$ is an $$(\varepsilon_1,\ldots,\varepsilon_m)\text{-superficial sequence for } J, I_1,\ldots,I_s \text { as in [TV]}.$$ Then $x_1,\ldots,x_m$ is a superficial sequence  of $(I_1,\ldots, I_s)$ with respect to $J$ by Remark 4.3 (ii).  Hence  Theorem 4.5 covers the main  result of Trung and Verma [TV]. 

 The notion of superficial elements goes back to P. Samuel [ZS]. The classical theory of superficial elements is an important tool in local algebra and has been continually developed (see e.g. [K],[HS], [RV]). Recall that
   $x$ 
 is called an $I_i$-superficial element of $R$ with respect to $( I_1,\ldots,I_s)$ if $x \in I_i$ and there exists a non-negative integer $c$ such that
$$(I_1^{n_1}\cdots I_i^{n_i+1}\cdots I_s^{n_s}: x)\bigcap
I_1^{n_1}\cdots I_i^{c}\cdots I_s^{n_s} = I_1^{n_1}\cdots I_s^{n_s}$$ 
for all $n_i \ge c$ and for all non-negative integers $ n_1, \ldots, n_{i - 1}, n_{i + 1}, \ldots, n_s.$       

\vskip 0.2cm

\hspace*{-0.6cm}{\bf Remark 4.7.} Let $J=I_0$ be an $\frak n$-primary ideal.
 It can be verified that if $x$ 
 is an $I_i$-superficial element of $R$ with respect to $(J, I_1,\ldots,I_s)$ then $$(J^{n+1}I_1^{n_1}\cdots I_i^{n_i+1}\cdots I_s^{n_s}:x) \bigcap J^{n}I_1^{n_1}\cdots I_s^{n_s} = J^{n+1}I_1^{n_1}\cdots I_s^{n_s}$$ for $n, n_1,\ldots,n_s\gg0.$ Indeed, since $J$ is $\frak n$-primary, there exists $u > c$ such that $I_i^{n_i} \subseteq JI_i^{c}$ for all $n_i \ge u.$ Consequently, $J^{n}I_1^{n_1}\cdots I_s^{n_s} \subseteq J^{n+1}I_1^{n_1}\cdots I_i^{c}\cdots I_s^{n_s}$ for all $n_i \ge u.$ Since $$(J^n I_1^{n_1}\cdots I_i^{n_i+1}\cdots I_s^{n_s}: x)\bigcap
J^nI_1^{n_1}\cdots I_i^{c}\cdots I_s^{n_s} = J^nI_1^{n_1}\cdots I_s^{n_s},$$ it follows that 
$$\begin{array}{ll}
&(J^{n+1} I_1^{n_1}\cdots I_i^{n_i+1}\cdots I_s^{n_s}: x)\bigcap
J^nI_1^{n_1}\cdots I_i^{n_i}\cdots I_s^{n_s} \\&\subseteq             
(J^{n+1} I_1^{n_1}\cdots I_i^{n_i+1}\cdots I_s^{n_s}: x) \bigcap J^{n+1}I_1^{n_1}\cdots I_i^{c}\cdots I_s^{n_s}\\&= J^{n+1}I_1^{n_1}\cdots I_i^{n_i}\cdots I_s^{n_s}
\end{array}$$
for all $n_i \ge u.$ So $$(J^{n+1}I_1^{n_1}\cdots I_i^{n_i+1}\cdots I_s^{n_s}:x) \bigcap J^{n}I_1^{n_1}\cdots I_s^{n_s} = J^{n+1}I_1^{n_1}\cdots I_s^{n_s}$$ for all  $n, n_1,\ldots,n_s\gg0.$
Hence if $x$ is  an $I_i$-superficial element of $R$ with respect to $(J, I_1,\ldots,I_s)$ then $x$ satisfies the condition (i) of a superficial element of $(I_1,\ldots, I_s)$ with respect to $J.$  
 Now we need to prove that if $I_1,\ldots,I_s$ are $\frak n$-primary, then $x$ satisfies also the condition (ii) of a superficial element of $(I_1,\ldots, I_s)$ with respect to $J.$ Indeed, 
by [HS, Lemma 17.2.4],
$$J^n I_1^{n_1}\cdots I_s^{n_s}:x=(0:x)+J^n I_1^{n_1}\cdots I_i^{n_i-1}\cdots I_s^{n_s}$$ for all $n, n_1,\ldots,n_s\gg0.$
Multiplying the last equation by $x$ yields
$$J^n I_1^{n_1}\cdots I_s^{n_s}\bigcap (x) = xJ^n I_1^{n_1}\cdots I_i^{n_i-1}\cdots I_s^{n_s}$$
for $n, n_1,\ldots,n_s\gg0.$ Hence $x$ satisfies the condition (ii) of a superficial element of $(I_1,\ldots, I_s)$ with respect to $J.$   
So if $I_1,\ldots,I_s$ are $\frak n$-primary and $x$ 
 is  an $I_i$-superficial element, then $x$ is  a superficial element of $(I_1,\ldots, I_s)$ with respect to $J.$

\vskip 0.2cm

\hspace*{-0.6cm}{\bf Remark 4.8.}
Let $J=I_0, I_1,\ldots,I_s$ be $\frak n$-primary ideals and $x$ 
  an $I_i$-superficial element. Then  $$\dim R/(x): I^\infty = \dim R/(x)= \dim R-1$$ and for any $e(J^{[k_0+1]},I_1^{[k_1]},\ldots,I_s^{[k_s]};R)$ with $(k_0+k_1+\cdots+k_s = \dim R -1= d-1),$ there exist $x_1,\ldots,x_n,\ldots, x_{d}$  of $R$ consisting of $k_1$ elements of $I_1,\ldots,k_i$ elements of $I_i,\ldots,k_s$ elements of $I_s,$ $x_{n+1},\ldots, x_{d}$ are $k_0+1$ elements of $J$ such that sequence $x_1,\ldots,x_{d}$ is a superficial sequence. Hence by Proposition 4.4(iii), $$e(J^{[k_0+1]},I_1^{[k_1]},\ldots,I_s^{[k_s]};R)=e(\bar J^{[1]},\bar I_1^{[0]},\ldots,\bar I_s^{[0]};\ovl R),$$ where $\ovl R=R/(x_1,\ldots,x_{d-1}),\bar J=J\ovl R,\bar I_1=I_1\ovl R,\ldots,\bar I_s=I_s\ovl R.$ Since $\dim \ovl R=1,$ $$e(\bar J^{[1]},\bar I_1^{[0]},\ldots,\bar I_s^{[0]};\ovl R)=l_R\biggl(\dfrac{R}{(x_1,\ldots,x_d)}\biggl)
-l_R\biggl(\dfrac{(x_1,\ldots,x_{d-1}):x_{d}}{(x_1,\ldots,x_{d-1})}\biggl)$$
by [HS, Theorem 17.4.6].
Since $x_1,\ldots,x_{d}$ is an $\frak n$-filter-regular sequence,  
$$l_R\biggl(\dfrac{R}{(x_1,\ldots,x_d)}\biggl)
-l_R\biggl(\dfrac{(x_1,\ldots,x_{d-1}):x_{d}}{(x_1,\ldots,x_{d-1})}\biggl)=
e(x_1,\ldots,x_d;R)$$
by [AB]. So we get the result of Risler and Teissier in [Te].\\
                
\vskip 0.2cm
Finaly, from  the above content  we would like to give the following conclusion.     
\vskip 0.2cm
\hspace*{-0.6cm}{\bf Remark 4.9.} It should be noted that all results of this section still are true if the condition (ii) of Definition 4.2 is replaced  by the  condition:     
$$(x)\bigcap {I_0}^{n_0} \cdots I_i^{n_i}\cdots I_s^{n_s} 
\equiv x{I_0}^{n_0}\cdots I_i^{n_i-1}\cdots I_s^{n_s}(\mathrm{mod }\;{I_0}^{n_0+1} \cdots I_i^{n_i}\cdots I_s^{n_s})\:$$
for $n_0,n_1,\ldots,n_s\gg0$ (see Remark 4.1).\\\\

\vspace*{10pt}
 
\centerline{\large\bf References}

\vspace*{6pt}
{\small {\begin{itemize}
\item [{[AB]}] A. Auslander and D. A. Buchsbaum, {\it Codimension and multiplicity}, Ann. Math. 68(1958), 625-657.
\item [{[Ba]}] P. B. Bhattacharya, {\it The Hilbert function of two ideals}, Proc. Cambridge. Philos. Soc. 53(1957), 568-575.
\item [{[BS]}] M. Brodmann and R. Y. Sharp, {\it Local cohomology: an algebraic  introduction with  geometric applications}, Cambridge studies in advanced mathematics, No 60, Cambridge University Press 1998.
\item[{[HHRT]}] M. Herrmann, E.  Hyry,  J.  Ribbe, Z. Tang,  {\it Reduction numbers and multiplicities of multigraded structures},  J. Algebra 197(1997), 311-341.
\item[{[HS]}] C. Huneke and I. Swanson, {\it Integral Closure of Ideals, Rings, and Modules}, London Mathematical Lecture Note Series 336, Cambridge University Press (2006).
\item[{[Hy]}]E.  Hyry, {\it The diagonal subring and the Cohen-Macaulay property of a multigraded ring}, Trans. Amer. Math. Soc. 351(1999), 2213-2232.
\item [{[K]}] D. Kirby, {\it A note on superficial elements of an ideal in a local ring}, Quart. J. Math. Oxford (2), {\bf 14} (1963), 21-28.
\item [{[KV]}] D. Katz,  J. K. Verma, {\it Extended Rees algebras and mixed multiplicities}, Math. Z. 202(1989), 111-128.
\item [{[KR1]}] D. Kirby and D. Rees, {\it Multiplicities in graded rings I: the general theory}, Contemporary Mathematics 159(1994), 209 - 267. 
\item [{[KR2]}] D. Kirby and D. Rees, {\it Multiplicities in graded rings II: integral equivalence and the Buchsbaum - Rim multiplicity}, Math. Proc. Cambridge Phil. Soc. 119 (1996),  425 - 445. 
 \item [{[KT]}]  S. Kleiman and A. Thorup, {\it Mixed Buchsbaum - Rim multiplicities}, Amer. J. Math. 118(1996), 529-569.
\item [{[MV]}] N. T. Manh and D. Q. Viet, {\it Mixed  multiplicities of modules over Noetherian local rings}, Tokyo J. Math. Vol. 29. No. 2, (2006), 325-345. 
\item [{[NR]}] D. G.  Northcott,  D. Rees, {\it Reduction of ideals in local rings},  Proc. Cambridge Phil. Soc. 50(1954), 145 - 158. 
\item [{[Re]}] D. Rees, {\it Generalizations of reductions and mixed multiplicities}, J. London. Math. Soc. 29(1984), 397-414.
\item [{[Ro]}] P. Roberts, {\it Local Chern classes, multiplicities and perfect complexes}, Memoire Soc. Math. France 38(1989), 145 - 161.
\item [{[RV]}] M. E. Rossi and G. Valla, {\it Hilbert function of filtered modules}, arXiv: 0710.2346 (2008).
\item [{[SV]}] J. Stuckrad and  W. Vogel, {\it Buchsbaum  rings and applications}, VEB Deutscher Verlag der Wisssenschaften. Berlin, 1986. 

\item [{[Sw]}] I. Swanson, {\it Mixed multiplicities, joint reductions and quasi-unmixed local rings }, J. London Math. Soc. 48(1993), no. 1, 1 - 14.

\item [{[Te]}] B. Teisier, {\it Cycles \`evanescents, sections planes, et conditions de Whitney}, Singularities \`a  Carg\'ese, 1972. Ast\'erisque,
\item [{[Tr1]}] N. V. Trung, {\it Reduction exponents and degree bound for the defining equation of graded rings}, Proc. Amer. Mat. Soc. 101(1987), 229-234.
\item [{[Tr2]}] N. V. Trung, {\it Positivity of mixed multiplicities}, J. Math. Ann. 319(2001), 33 - 63.
\item [{[TV]}] N. V. Trung and  J. Verma, {\it   Mixed  multiplicities of ideals versus mixed volumes of polytopes}, Trans. Amer. Math. Soc. 359(2007), 4711-4727.
\item [{[Ve]}] J. K. Verma, {\it Multigraded  Rees algebras and mixed multiplicities}, J. Pure and  Appl.  Algebra 77(1992), 219-228.
\item [{[Vi]}] D. Q. Viet, {\it Mixed multiplicities of arbitrary ideals in local rings}, Comm. Algebra. 28(8)\\(2000), 3803-3821.

\item [{[ZS]}] O. Zariski and P. Samuel, {\it Commutative Algebra}, Vol II, Van Nostrand, New York, 1960.

\end{itemize}}

\end{document}